\documentclass{article}
\usepackage{amsfonts}
\usepackage{amsmath}
\usepackage{a4wide}
\usepackage{amssymb}

\newlength{\tewi}
\newlength{\tehe}
\newtheorem{lem}{Lemma}[section]
\newtheorem{theor}[lem]{Theorem}
\newtheorem{corol}[lem]{Corollary}
\newtheorem{propo}[lem]{Proposition}
\newtheorem{remar}[lem]{Remark}
\newtheorem{defin}[lem]{Definition}
\newtheorem{conje}[lem]{Conjecture}

\newenvironment{theorem}
{\begin{theor}\sl} {\end{theor}}

\newenvironment{rtheorem}[1]
{\begin{theor}[{\rm #1}]\sl} {\end{theor}}

\newenvironment{proposition}
{\begin{propo}\sl } {\end{propo}}

\newenvironment{definition}
{\begin{defin}\rm } {\end{defin}}

\newcommand\proof{\noindent{\bf Proof}\quad}

\newcommand\qed{\hfill$\blacksquare$}

\newcommand{\bau}[1]{{\sc #1}}
\newcommand{\bjo}[1]{{\it #1}}
\newcommand{\bvo}[1]{{\bf #1}}
\newcommand{\bbo}[1]{{\it #1}}
\newcommand{\bpa}[1]{{#1}}

\newcommand\barK{{\bar K}}
\newcommand\cala{{\cal A}}
\newcommand\calf{{\cal F}}

\newcommand\cals{{\cal S}}

\newcommand\gal{{\rm Gal}}

\newcommand\gera{{\mathfrak a}}

\newcommand\norm{{\cal N}}

\newcommand\ord{{\rm Ord}}

\newcommand\Q{{\mathbb Q}}
\newcommand\R{{\mathbb R}}
\newcommand\Z{{\mathbb Z}}

\newcommand\eps\varepsilon

1
\newcommand\addtoc{\addtocontents{toc}{\protect\vspace{-0.8\baselineskip}}}

\title{\vskip-2\baselineskip
Divisibility of class numbers: enumerative approach}
\author{Yuri F. Bilu and Florian Luca}

\begin{document}

\hfuzz 2pt

{\sloppy

\maketitle

{\footnotesize \tableofcontents }

\section{Introduction}
\addtoc It is well-known since Gauss that infinitely many
quadratic fields have even class number. In fact, if~$K$ is a
quadratic field of discriminant~$D$, having~$r$ prime divisors,
then the class number~$h_K$ is divisible by~$2^{r-1}$ if ${D<0}$
and by~$2^{r-2}$ if ${D>0}$. See~\cite[Theorem~3.8.8]{BS66} for a
more precise statement.

In~1922 Nagell \cite[Satz~VI]{Na22} obtained the following
remarkable result: {\sl given a positive integer~$\ell$, there
exist infinitely many imaginary quadratic fields with class number
divisible by~$\ell$}. See~\cite{AC55} for a different proof.

It took almost half a century to extend Nagell's result to real
quadratic field, see Yamamoto~\cite{Ya70} and
Weinberger~\cite{We73}. Uchida~\cite{Uc74} extended this to cyclic
cubic fields. In mid-eighties, Azuhata and Ichimura~\cite{AI84}
and Nakano~\cite{Na83,Na85} obtained similar results for fields of
arbitrary degree~$n$.

Recently Murty~\cite{Mu99} gave  quantitative versions of the
theorems of Nagell and Yamamoto-Weinberger on quadratic fields. He
proved that for all sufficiently large~$X$ there exist at
least~${c(\ell)X^{1/2+1/\ell}}$ imaginary quadratic fields and at
least ${c(\ell,\eps) X^{1/{4\ell}-\eps}}$ real quadratic fields
with class number divisible by~$\ell$ and discriminant not
exceeding $X$ in absolute value.  (The second exponent can be
replaced by~${{1/{2\ell}-\eps}}$ if~$\ell$ is odd.) Various
refinement and extensions of Murty's results were suggested in
\cite{CM03,Lu03,So00,Yu02}.

Much less is known for fields of higher degree. In~\cite{HL03}, it
is shown that at least ${c(\ell)X^{1/{6\ell}}/\log X}$ pure cubic
fields have discriminant not exceeding~$X$ and class number
divisible by~$\ell$.

In this paper, we extend these results to fields of  degree ${n\ge
3}$.

\begin{theorem}
\label{tm} Let~$n$ and~$\ell$ be positive integers, ${n\ge 3}$,
and put ${\mu = \frac1{2(n-1)\ell}}$. There exist positive real
numbers ${X_0=X_0(n, \ell)}$  and ${c=c(n,\ell)}$ with the
following property. For any ${X>X_0}$ there is at least $cX^\mu$
pairwise non-isomorphic number fields of degree~$n$, discriminant
not exceeding~$X$, and the class number divisible by~$\ell$.
\end{theorem}

Actually, we prove slightly more: for all those fields the class
group has an element of exact order~$\ell$.

The famous {\it Cohen-Lenstra heuristics}~\cite{CL84,CM90} predict
that number fields of degree ${n>1}$ with class number divisible
by~$\ell$ should have positive density among all number fields of
degree~$n$. More precisely, denote by $\calf_n (X)$ the set of all
non-isomorphic number fields of degree~$n$ and discriminant not
exceeding~$X$ and put ${\calf_n^{(\ell)}(X) =\{K \in \calf_n(X):
\ell|h(K)\}}$. Then, as ${X\rightarrow\infty}$, the quotient
${\left |\calf_n^{(\ell)}(X)\right|/\left |\calf_n(X)\right|}$
(conjecturally) tends to a positive rational number, which can be
explicitly expressed in terms of certain finite Euler-type
products.

This conjecture seems to be out of reach at the present state of
our knowledge. Theorem~\ref{tm} implies that number fields of
degree~$n$ with class group divisible by~$\ell$ have positive {\it
logarithmic} lower density among all number fields of degree~$n$~:
\begin{equation}
\label{equot} \liminf_{X\rightarrow \infty} \frac{\log\left
|\calf_n^{(\ell)}(X)\right|}{\log\left |\calf_n(X)\right|}\ge
\frac2{(n-1)(n+2)\ell}.
\end{equation}
This is an immediate consequence of Theorem~\ref{tm} and the
inequality
\begin{equation}
\label{edib} \left |\calf_n(X)\right| \le c(n) X^{(n+2)/4}
\end{equation}
due to Schmidt~\cite{Sc95} (see also
\cite[Proposition~9.3.4]{Co00}). For large~$n$
inequality~(\ref{equot}) can be refined due to the recent work of
Ellenberg and Venkatesh~\cite{EV04}.

\medskip

The argument of the present paper relies on the famous
construction of Ankeny-Brauer-Chowla fields~\cite{ABC56} and is
strongly inspired by the work of Sprindzhuk
\cite[Section~8.6]{Sp82} and Halter-Koch {\it et
al.}~\cite{HLPT99}. In Sections~\ref{sthin}--\ref{sabc} we collect
necessary facts about thin sets and Ankeny-Brauer-Chowla fields.
The proof of Theorem~\ref{tm} occupies
Sections~\ref{smp}--\ref{sabcsi}.

\subsection{Notations and Conventions} {\bf All fields in this
paper are of characteristic~$0$.} Let~$K$ be a field. We
write~$\barK$ for
its algebraic closure. By the Galois group
of a separable polynomial ${f(x) \in K[x]}$ we mean the Galois
group of the splitting field of~$f$ over~$K$, realized as a
subgroup of the symmetric group~$\cals_n$, where ${n=\deg f}$. In
particular,~$f$ is irreducible over~$K$ if and only if its Galois
group is transitive.

Unless the contrary is stated explicitly (as it is done in
Section~\ref{sasts}), small letters~$t,x,y,z$ with or without
indices denote indeterminates algebraically independent over the
base field.

\paragraph{Acknowledgements}
We thank Henri Cohen, Jacques  Martinet and  Michel Olivier for
useful discussions. We are indebted to W\l adys\l aw Narkiewicz
and the referee for drawing our attention to the work of Nagell
and Nakano.  Yuri Bilu thanks Vera Bergelson for inspiring
comments. This work was supported in part by the Joint Project
France-Mexico ECOS--ANUIES M02--M01.

\section{Thin sets}
\label{sthin} \addtoc Let~$K$ be a field and let~$n$ be a positive
integer. Let~$\Upsilon$ be a subset of the affine space~$K^n$. The
set~$\Upsilon$ is called {\it basic thin set of the first type} if
there exists a {\bf non-zero} polynomial ${F(\underline t) \in
K(\underline t})$  (where ${\underline t= (t_1,\ldots,t_n)}$) such
that ${(\underline \tau) \in \Upsilon}$ if and only if
${F(\underline \tau)=0}$. The set~$\Upsilon$ is a {\it basic thin
set of the second type} if there exists an $K$-irreducible
polynomial ${F(\underline t, x)\in K(\underline t, x)}$ with
${\deg_xF\ge 2}$ such that ${(\underline \tau) \in \Upsilon}$ if
and only if the specialized polynomial ${F(\underline \tau, x)\in
K[x]}$ has a root in~$K$. The set~$\Upsilon$ is called {\it thin}
if it is contained in a finite union of basic thin sets. It is
obvious that the union of finitely many thin sets is thin, and
that a subset of a thin set is thin.

Serre \cite[Section~9.1]{Se97} gives a differently looking, but
equivalent definition of thin sets.

This following property must be known, but we could not find it in
the literature.

\begin{proposition}
\label{pinter} Let~$L$ be a finitely generated extension of the
field~$K$, and~$\Upsilon$ a thin subset of~$L^n$. Then ${\Upsilon
\cap K^n}$ is a thin subset of~$K^n$.
\end{proposition}

\proof The case of finite extension~$L/K$ is considered in
\cite[page~128]{Se97}, so we are left with the pure transcendent
case. Thus, assume that ${L=K(\underline z)}$, where ${\underline
z=(z_1, \ldots, z_s)}$, and let ${\Upsilon\subset L^n}$ be a basic
thin set of the first type, defined by the polynomial
${F(\underline t)\in L[\underline t]}$.

We may assume that~$K^{n+s}$ is not a thin subset of itself;
otherwise~$K^n$ is a thin subset of itself as well (cf.
\cite[Section~9.4]{Se97}), and the statement becomes trivial. It
follows that~$F$, viewed as a rational function in ${\underline t,
\underline z }$, is defined and does not vanish at certain
${\left(\underline \tau', \underline \zeta\right)\in K^{n+s}}$.
Hence, the polynomial ${F_{\underline\zeta}(\underline t) \in
K[\underline t]}$, obtained from~$F$ by specialization
${\underline z=\underline\zeta}$, is defined and non-zero.  For
any ${\underline \tau \in \Upsilon\cap K^n}$ we have
${F_{\underline\zeta}(\underline \tau)=0}$. Hence, ${\Upsilon \cap
K^n}$ is thin.

One argues similarly in the case when ${\Upsilon\subset L^n}$ is a
basic thin set of the second type, defined by the polynomial
${F(\underline t, x)\in L[\underline t, x]}$. This time, we find
${\underline\zeta \in K^s}$ such that  the polynomial
${F_{\underline\zeta}(\underline t, x)\in K[\underline t, x]}$ has
no factors of $x$-degree~$1$. Let ${F_{\underline \zeta}=G_1\ldots
G_k}$ be the irreducible decomposition of~$F_{\underline\zeta}$ in
$K[\underline t, x]$. Then every~$G_i$ is of $x$-degree at
least~$2$, and ${\Upsilon \cap K^n}$ lies in the union of the
corresponding basic thin sets of the second type.  \qed

\begin{theorem}
\label{tgal} \label{tirr} Let   ${F(t_1, \ldots, t_n, x)\in
K[t_1,\ldots,t_n, x]}$ be a polynomial of $x$-degree~$m$, and let
${s \le n}$. Let ${G\le \cals_m}$ be the Galois group of~$F$ over
the field ${K(t_1,\ldots,t_n)}$. Then for all ${(\tau_1, \ldots ,
\tau_s) \in K^s}$ outside a thin set the polynomial ${F(\tau_1,
\ldots, \tau_s, t_{s+1}, \ldots, t_n, x)\in K[t_{s+1}, \ldots,
t_n, x]}$ is separable, of $x$-degree~$m$, and its Galois group
over ${K(t_{s+1}, \ldots, t_n)}$ is~$G$.

In particular, if~$F$ is irreducible over ${K(t_1,\ldots,t_n)}$,
then for all ${(\tau_1, \ldots , \tau_s) \in K^s}$ outside a thin
set the polynomial ${F(\tau_1, \ldots, \tau_s, t_{s+1}, \ldots,
t_n, x)\in K[t_{s+1}, \ldots, t_n, x]}$ is irreducible over
${K(t_{s+1}, \ldots, t_n)}$.
\end{theorem}

\proof The case ${s=n}$ is treated  in \cite[Section~9.2,
Propositions~1 and~2]{Se97}. The general case reduces to ${s=n}$
by Proposition~\ref{pinter}.\qed

\medskip

Hilbert's irreducibility theorem asserts that~$K^n$ is not thin
for a finitely generated field~$K$. We shall use the following
quantitative version for ${K=\Q}$, due to S.~Cohen~\cite{Co81}.
See also Serre \cite[Section~13.1, Theorem~1]{Se97}.
\begin{theorem}
\label{thit} Let~$\Upsilon$ be a thin subset of~$\Q^n$. Then there
exists a positive constant ${c=c(\Upsilon)}$ such that  for
${X>1}$ we have
$$
\left|\Upsilon \cap \Z^n\cap [-X,X]^n\right|\le cX^{n-1/2}\log X.
\eqno \blacksquare
$$
\end{theorem}
(For ${n=1}$ the $\log$-factor can be omitted.)

\section{A special thin set}
\addtoc \label{sasts} In this (and only this) section we use
capital letters ${X,Y,Z,\ldots}$ for independent variables,
reserving small letters $x,y,z\ldots$ for algebraic functions.

\begin{proposition}
\label{puv}
Let~$K(x)$ be the field of rational functions over~$K$. Let
${u,v\in K(x)}$ satisfy the following:~$v$ has a simple zero (or a
simple pole) in~$\barK$ which is neither a zero nor a pole of~$u$.
Then for any positive integers~$k$ and~$\ell$ we have ${\left[K(x,
u^{1/k}, v^{1/\ell})\colon K(x, u^{1/k})\right]=\ell}$.
\end{proposition}

\proof Obviously, ${\left[K(x, u^{1/k}, v^{1/\ell})\colon K(x,
u^{1/k})\right]\le\ell}$, so it remains to prove that
$$
\left[K(x, u^{1/k}, v^{1/\ell})\colon K(x, u^{1/k})\right]\ge\ell.
$$
Let ${\alpha\in \barK}$ be a simple zero (or pole) of~$v$, which
is neither a zero nor a pole of~$u$, and let $\ord_\alpha$ be the
corresponding  place of the field~$K(x)$. This place is unramified
in the field ${K(x, u^{1/k})}$ , but it is ramified in the field
${K(x, u^{1/k}, v^{1/\ell})}$, with ramification index~$\ell$.
Hence ${\left[K(x, u^{1/k}, v^{1/\ell})\colon K(x,
u^{1/k})\right]\ge\ell}$, as wanted.\qed

\medskip

\begin{proposition}
\label{pkum} Let~$K$ be a field. Consider the polynomial
$$
f(T,X) := \left(X-a_1\right)\cdots \left(X-a_{n-1}\right) \left
(X-\alpha T^\ell-\beta\right)-1 \in K[X,T],
$$
where~$n, \ell$ are positive integers and ${a_1,\ldots, a_{n-1},
\alpha, \beta \in K^\ast}$.
Let ${\nu_1, \ldots, \nu_{n-1}, \nu}$ be integers, ${\nu \ne0}$.
Let~$\Upsilon$ be the subset of~$K$ defined as follows: ${\tau\in
K}$ belongs to~$\Upsilon$ if for some root~$\xi$ of the polynomial
${f(\tau, X) \in K[X]}$, and for some determination of ${\zeta =
\left(\xi(\xi-a_1)^{\nu_1/\nu}\cdots
(\xi-a_{n-1})^{\nu_{n-1}/\nu}\right)^{1/\ell}}$, we have ${[K(\xi,
\zeta)\colon K(\xi)]<\ell}$.

Assume that the polynomial ${f(T,X)}$ is irreducible over~$K$, and
that the polynomial ${f(0, X) \in K[X]}$ has a simple root
in~$\barK^\ast$. Then~$\Upsilon$ is thin.
\end{proposition}

\proof Let~$K(x)$ be the field of rational functions and let
${t,z\in \overline{K(x)}}$ be defined by
$$
t= \left(\frac{f(0,x)}{\alpha(x-a_1)\cdots
(x-a_{n-1})}\right)^{1/\ell}, \quad
z=\left(x(x-a_1)^{\nu_1/\nu}\cdots
(x-a_{n-1})^{\nu_{n-1}/\nu}\right)^{1/\ell}.
$$
Then ${f(t,x)=0}$, and, since~$f$ is irreducible, we have
\begin{equation}
\label{enell} [K(x,t)\colon K(t)]=n, \quad [K(x,t)\colon
K(x)]=\ell.
\end{equation}
Also,
\begin{equation}
\label{egell} [K(x, z)\colon K(x)]\ge \ell
\end{equation}
(consider the ramification at~$0$).

Further, by the assumption, ${f(0, x)}$ has a non-zero simple
root. This root is distinct from any of the numbers ${0, a_1,
\ldots, a_{n-1}}$. Proposition~\ref{puv} implies that
${[K(x,z,t)\colon K(x,z)]=\ell}$. Combining this
with~(\ref{enell}) and~(\ref{egell}), we obtain
${m:=[K(x,z,t)\colon K(t)]\ge n\ell}$.

Now let ${y\in K(x,z,t)}$ be such that ${K(x,z,t)=K(t,y)}$, and
let ${g(T,Y) \in K[T,Y]}$ be the irreducible polynomial over~$K$
such that ${g(t,y)=0}$. Then ${\deg_y g=m}$. Applying
Theorem~\ref{tirr}, we find  a thin set ${\Upsilon_1\subset K}$
such that for any ${\tau \in K\setminus \Upsilon_1}$, the
polynomial ${g(\tau, Y)\in K[Y]}$ is of degree~$m$ and irreducible
over~$K$.

On the other hand, there exists ${h(T,X,Z)\in K(T)[X,Z]}$ such
that ${y=h(t,x,z)}$. Denote by~$d(T)$ the denominator of
$h(T,X,Z)$.

Fix ${\tau\in \Upsilon}$, together with the corresponding~$\xi$
and~$\zeta$. Then ${[K(\xi,\zeta)\colon K]<n\ell\le m}$. Assume
that ${d(\tau)\ne 0}$. Then ${h(\tau, \xi, \zeta)}$ is a root of
$g(\tau, Y)$ of degree~${<m}$ over~$K$. Hence, either ${\deg
g(\tau, Y) < m}$ or $g(\tau, Y)$ is reducible over~$K$. In both
cases ${\tau \in \Upsilon_1}$.

We have proved that ${\Upsilon \subseteq \Upsilon_1 \cup
\{\mbox{the roots of $d(T)$}\}}$. Hence,~$\Upsilon$  is thin. \qed

\section{The Ankeny-Brauer-Chowla fields}
\label{sabc} \addtoc Let ${a_1, \ldots, a_n}$, where ${n\ge 3}$,
be pairwise distinct integers and ${f(x) = (x-a_1)\cdots
(x-a_n)-1}$. It is well-known that~$f(x)$ is an irreducible
polynomial \cite[Problem~8.121]{PS64}. The number fields, defined
by  such polynomials, are called {\it Ankeny-Brauer-Chowla
fields}~\cite{ABC56} (ABC-fields in the sequel).

Let~$\xi$ be a root of~$f$. The main property of the ABC-fields is
that, under mild assumptions about the numbers ${a_1, \ldots,
a_n}$, the field ${K=\Q(\xi)}$ is totally real, and the numbers
${\xi - a_1, \ldots, \xi-a_{n-1}}$ form a full rank system of
units of~$K$.

Below we summarize the properties of the Ankeny-Brauer-Chowla
polynomials and fields, to be used in this paper. In the sequel,
${a_1,\ldots, a_{n-1}}$ are fixed pairwise distinct integers,~$a$
runs in the set of integers distinct from any of ${a_1,\ldots,
a_{n-1}}$, and ${f_a(x) = (x-a_1)\cdots (x-a_{n-1})(x-a)-1}$.
Unless the contrary is stated explicitly, implicit constants in
this section depend only on ${a_1,\ldots, a_{n-1}}$. In
particular, {\it sufficiently large} means {\it exceeding a
positive constant depending on ${a_1,\ldots, a_{n-1}}$}.

\begin{proposition}
\label{pabc}
Assume that~$|a|$ is sufficiently large. Then we have
the following.

\begin{enumerate}
\item
\label{iroots} The polynomial~$f_a(x)$ has~$n$ real roots ${\xi_1,
\ldots, \xi_{n-1}, \xi}$ satisfying
\begin{align}
\label{eroots}
|\xi_k-a_k|&\ll |a|^{-1} \quad (k=1, \ldots, n-1),\\
|\xi -a|&\ll|a|^{1-n}.
\end{align}
In particular, the number field ${K_a:=\Q(\xi)}$ is totally real.

\item
\label{idis} The discriminant of the field~$K_a$
is~$O\left(|a|^{2(n-1)}\right)$.
\item
\label{iunits} The numbers ${\xi - a_1, \ldots, \xi- a_{n-1}}$, form a
full rank system of independent units of the field~$K_a$.

{\rm (These numbers are called {\it basic ABC-units}. The
multiplicative group, generated by the basic ABC-units, is called
{\it the group of ABC-units}.)}

\item
\label{iindex} Assume that the field~$K_a$ is
primitive\footnote{that is, it has no proper subfield distinct
from~$\Q$.}, and that the absolute value of its discriminant
exceeds~$|a|^\kappa$, where~$\kappa$ is a positive number. Then
the group of ABC-units is of index at most
$O\left(\kappa^{1-n}\right)$ in the group of all units.
\end{enumerate}
\end{proposition}
\proof Parts~\ref{iroots} and~\ref{idis} are obvious. To prove
Part~\ref{iunits}, consider the real embeddings
\begin{equation}
\label{eemb}
\begin{array}{r@{\,}c@{\,}l}
\sigma_i:K_a&\rightarrow&\R\\
\xi&\mapsto&\xi_i
\end{array}
\quad (i=1, \ldots, n-1).
\end{equation}
Then~(\ref{eroots}) implies that
${\log\left|\sigma_i(\xi-a_j)\right|=-\delta_{ij}\log|a|+O(1)}$,
where~$\delta_{ij}$ is the Kronecker symbol. Hence,
$$
R_{ABC}:=\left|\det\left[\log\left|\sigma_i(\xi-a_j)\right|\right]_{1\le
i,j\le n-1}\right|= (\log |a|)^{n-1} + O\left((\log
|a|)^{n-3}\right).
$$
In particular, ${R_{ABC}\ne 0}$ for sufficiently large~$|a|$,
which proves Part~\ref{iunits}.

For Part~\ref{iindex}, recall (cf.~\cite{Si84,Fr89}) that the
regulator~$R$ and the discriminant~$D$ of a primitive field~$K$
satisfy the inequality ${R\gg(\log |D|)^r}$, where~$r$ is the rank
of the unit group of~$K$ and the implicit constant depends on the
degree of~$K$. For the totally real field~$K_a$ we have ${r=n-1}$,
which, together with the assumption ${|D|\ge |a|^\kappa}$, imply
${R\gg(\kappa\log|a|)^{n-1}}$. Hence, ${R_{ABC}/R\ll
\kappa^{1-n}}$, as wanted. \qed

\medskip
(It might be pointed out that for sufficiently large~$|a|$ the
implicit constant in Part~\ref{iindex} depends only on~$n$. For
instance, using Theorem~C of Friedman~\cite{Fr89}, one can show
that for sufficiently large~$|a|$ the index of the ABC-units in the
group of all units does
not exceed ${Cn^{2n}\kappa^{1-n}}$, where~$C$ is an absolute
constant. We shall not use this more precise estimate.)

\medskip

Sprindzhuk~\cite[Lemma~8.6.4]{Sp82} showed that, for ${n \ge 3}$,
distinct ABC-fields are seldom isomorphic. Below, we reproduce his
result in a slightly refined form.

\begin{proposition} (Sprindzhuk)\quad
\label{pspr} Assume that ${n\ge 3}$. Let~$A$ be a sufficiently
large positive integer, and let~$S$ be a set of integers~$a$
satisfying ${A\le |a|\le 2A}$ and such that for all ${a \in S}$
the fields~$K_a$ are isomorphic to the same field~$K$. Then
${|S|\le n(n-1)(n-2)}$.
\end{proposition}

\proof Assume that  ${|S|> n(n-1)(n-2)}$. Since~$\R$ has
exactly~$n$ distinct subfields isomorphic to~$K$, the set~$S$ has
more than ${(n-1)(n-2)}$ elements~$a$ such that all the
fields~$K_a$ are the same. Further, let ${\sigma_i\colon
K_a\rightarrow \R}$ be defined as in~(\ref{eemb}). Then, for a
given~$K_a$, there exist ${(n-1)(n-2)}$ possibilities for the pair
${(\sigma_1, \sigma_2)}$. Hence, there are distinct~$a$ and~$a'$
such that ${K_a=K_{a'}}$, \  ${\sigma_1=\sigma_1'}$ and
${\sigma_2=\sigma_2'}$. (Here and below ${\xi', \xi_1',\ldots,
\xi_{n-1}', \sigma_1', \ldots, \sigma_{n-1}'}$ have the same
meaning for~$a'$ as ${\xi, \xi_1,\ldots, \xi_{n-1}, \sigma_1,
\ldots, \sigma_{n-1}}$ for~$a$.) It follows that ${\xi-\xi'}$ is a
non-zero algebraic integer from the field~$K_a$, and
${\sigma_i(\xi-\xi')=\xi_i-\xi_i'}$ for ${i=1,2}$.
Using~(\ref{eroots}) and the assumption ${A\le |a|, |a'|\le 2A}$,
we obtain ${|\xi-\xi'|\ll A}$, as well as
${|\sigma_i(\xi-\xi')|\ll A^{-1}}$ for ${i=1,2}$ and
${|\sigma_i(\xi-\xi')|\ll1}$ for ${i=3,\ldots,n-1}$. Hence
$$
1\le \left|\norm_{K_a}(\xi-\xi') \right|=
\left|\xi-\xi'\right|\prod_{i=1}^{n-1}|\sigma_i(\xi-\xi')|\ll
A^{-1},
$$
which is a contradiction for sufficiently large values of~$A$. \qed

\section{Construction of the main polynomial}
\label{smp} \addtoc Starting from section, we begin the proof of
Theorem~\ref{tm}. Until the end of the paper, we fix positive
integers~$n$ and~$\ell$. Unless the contrary is stated explicitly,
{\bf we shall always assume that ${n \ge 3}$}.

In this section, we construct, for the given~$n$ and~$\ell$, a
special polynomial in two variables, which will be used in the
subsequent sections to produce Ankeny-Brauer-Chowla fields having
required properties.
\begin{theorem}
\label{tmp}  There exists pairwise distinct non-zero integers
${a_1,\ldots,a_{n-1}}$ such that
the polynomial
\begin{equation}
\label{emp}
f(t,x):=\left(x-a_1\right)\cdots\left(x-a_{n-1}\right)\left(x-
(-1)^{n-1}\frac{t^\ell-1}{a_1\cdots a_{n-1}}\right) -1\in \Q[t,x]
\end{equation}
has a symmetric Galois group over the field~$\Q(t)$, and the
polynomial ${f(0,x)}$ is separable.
\end{theorem}

Our starting point is the following result of Halter-Koch {\it et
al.} \cite[Proposition~3.1]{HLPT99}:

\begin{proposition}
\label{phlpt} Let~$K$ be a field, ${\gamma \in K^\ast}$ and
${t_1,\ldots,t_n}$ (algebraically independent) indeterminates
over~$K$. Then the Galois group of the polynomial  ${(x-t_1)\cdots
(x-t_n)-\gamma}$ over $K(t_1,\ldots,t_n)$ is~$\cals_n$. \qed
\end{proposition}

\begin{proposition}
\label{pgal} Let~$F$ be a field and~$H$  a finite Galois extension
of~$F$ with Galois group~$\cals_n$, where ${n\ge 4}$.
Let~$\alpha$ an element of~$\bar F$ such that ${\alpha^\ell\in
F}$. Then ${\gal(H(\alpha)/F(\alpha))}$ is either~$\cals_n$ or the
alternating group $\cala_n$.
\end{proposition}

\proof Let~$\zeta$ be a primitive $\ell$-th root of unity and put
${F_1:=F(\alpha, \zeta)}$ and ${H_1:=H(\alpha, \zeta)}$. Since
$$
\gal(H_1/F_1) \le \gal(H(\alpha)/F(\alpha))\le \gal(H/F)=\cals_n,
$$
it suffices to show that ${\gal(H_1/F_1)\ge \cala_n}$.

Since both~$H_1$ and~$F_1$ are Galois extensions of~$F$, the group
 ${\gal(H_1/F_1)= \gal \left(H/(H\cap
F_1)\right)}$ is an invariant subgroup of ${\cals_n= \gal(H/F)}$.
And it cannot be trivial because $\gal(F_1/F)$ is a meta-abelian
group, while~$\cals_n$ for ${n\ge 4}$ is not. It follows that
${\gal(H_1/F_1)\ge \cala_n}$, as wanted.\qed

\begin{proposition}
\label{psn}
Let~$K$ be a field, ${\alpha, \gamma \in K^\ast}$,
${\beta \in K}$ and ${n\ge 4}$. Then the Galois group of the
polynomial
\begin{equation}
\label{epoly}
\left(x-t_1\right)\cdots\left(x-t_{n-1}\right)\left(x-
\frac{\alpha t^\ell+\beta}{t_1\cdots t_{n-1}}\right) -\gamma
\end{equation}
over $K(t_1,\ldots,t_{n-1}, t)$ is~$\cals_n$.
\end{proposition}
\proof Put
$$
t_n:= \frac{\alpha t^\ell+\beta}{t_1\cdots t_{n-1}}.
$$
Proposition~\ref{phlpt} implies that the Galois group of
polynomial~(\ref{epoly}) over $K(t_1,\ldots, t_n)$ is~$\cals_n$.
Using Proposition~\ref{pgal}, we conclude that the Galois group
of~(\ref{epoly}) over $K(t_1,\ldots,t_{n-1}, t)$ is~$\cals_n$
or~$\cala_n$.

It remains to show that the $x$-discriminant of~(\ref{epoly}) is
not a square in $K(t_1,\ldots,t_{n-1}, t)$. It suffices to verify
that the $x$-discriminant~$D(t)$ of the polynomial
${g(t,x)=(x-1)^{n-1}(x-\alpha t^\ell-\beta) -\gamma}$ (obtained
from~(\ref{epoly}) by specializing ${t_1=\ldots= t_{n-1}=1}$) is
not a square in~$K(t)$.

Put ${a(t) = \alpha t^\ell+\beta-1}$, so that ${g(t)
=(x-1)^{n-1}(x-1-a(t))-\gamma}$. Then
$$
\frac {\partial g}{\partial x} (x,t) =
n(x-1)^{n-2}\left(x-1-\frac{n-1}n a(t)\right),
$$
whence
$$
D(t) =n^ng(t,1)^{n-2}g\left(t, \frac{n-1}n a(t)+1\right)=
(-1)^{n-1}\gamma^{n-2}\left((n-1)^{n-1}a(t)^n+n^n\gamma\right).
$$
Thus, ${\deg D(t)=n\ell}$ and ${D'(t) = \delta
a(t)^{n-2}t^{\ell-1}}$, where ${\delta \in K^\ast}$. Since~$D(t)$
does not vanish at  the roots of~$a(t)$, the only possible
multiple root of~$D(t)$ is~$0$, and if it is, its multiplicity
is~$\ell$. Hence,~$D(t)$ is not a square of a polynomial, as
wanted.\qed

\begin{proposition}
\label{psep} Let~$K$ be a field, ${\beta \in K}$ and ${\gamma \in
K^\ast}$. Assume that\footnote{This assumption can be dropped, but
the argument would become slightly more involved.}
\begin{equation}
\label{enn}
(n-1)^{n-1}(\beta-1)^n + n^n\gamma \ne 0
\end{equation}
Then the
polynomial
\begin{equation}
\label{epoly0}
\left(x-t_1\right)\cdots\left(x-t_{n-1}\right)\left(x-
\beta{t_1^{-1}\cdots t_{n-1}^{-1}}\right) -\gamma
\end{equation}
is separable over $K(t_1,\ldots,t_{n-1})$.
\end{proposition}
\proof Again, it suffices to show that the polynomial
${g(x)=(x-1)^{n-1}(x-\beta)-\gamma}$ (obtained from~(\ref{epoly0})
by specializing ${t_1=\ldots= t_{n-1}=1}$) is separable over~$K$.
We have
$$
g'(x)= n(x-1)^{n-2}\left(x-\frac{(n-1)\beta +1}n\right).
$$
Since ${g(1)=-\gamma\ne 0}$ and ${g\left(\frac{(n-1)\beta
+1}n\right)\ne 0}$ by~(\ref{enn}), the result follows. \qed

\medskip\noindent
{\bf Proof of Theorem~\ref{tmp}}\quad One immediately verifies
that $f_3(t,x)={(x^2-1)(x+t^\ell-1)-1}$ is an irreducible
over~$\Q(t)$ polynomial in~$x$,  and its $x$-discriminant is not a
square in~$\Q(t)$. Hence, its Galois groups over~$\Q(t)$
is~$\cals_3$. Since  ${f_3(0,x)}$ is separable, this proves the
theorem for ${n=3}$.

Assume now that ${n\ge 4}$ and consider the polynomial
$$
F\left (t_1, \ldots, t_{n-1}, t, x \right)=
\left(x-t_1\right)\cdots\left(x-t_{n-1}\right)\left(x-
(-1)^{n-1}\frac{ t^\ell-1}{t_1\cdots t_{n-1}}\right) -1 \in
\Q\left (t_1, \ldots, t_{n-1}\right)[ t, x].
$$
Propositions~\ref{psn} and~\ref{psep} imply that the Galois group
of~$F$ over $\Q(t_1, \ldots, t_{n-1}, t)$ is~$\cals_n$, and that
${F\left (t_1, \ldots, t_{n-1}, 0, x \right)}$ is separable over
$\Q(t_1, \ldots, t_{n-1})$.

By Theorem~\ref{tgal}, there exists a thin set ${\Upsilon
\subseteq \Q^{n-1}}$ such that for any ${(\tau_1, \ldots,
\tau_{n-1})}\in {(\Q^\ast)^{n-1}\setminus \Upsilon}$, the Galois
group of the specialized polynomial $F\left (\tau_1, \ldots,
\tau_{n-1}, t, x \right)$ is~$\cals_n$, and the polynomial $F\left
(\tau_1, \ldots, \tau_{n-1}, 0, x \right)$ is separable. Finally,
Theorem~\ref{thit} implies that there exist pairwise distinct
non-zero integers ${a_1, \ldots, a_{n-1}}$ such that ${(a_1,
\ldots, a_{n-1})\notin\Upsilon}$. This completes the proof of the
theorem.\qed

\section{Suitable integers}
\addtoc Recall that we fix positive integers~$n$ and~$\ell$ with
${n \ge 3}$. In addition, starting from this point, we fix, once
and for all, pairwise distinct non-zero integers ${a_1,\ldots,
a_{n-1}}$ (which exist by Theorem~\ref{tmp}) such that the
polynomial ${f(t,x)}$, defined in~(\ref{emp}), has Galois
group~$\cals_n$ over~$\Q(t)$, and the polynomial ${f(0,x)}$ is
separable.  Unless the contrary is stated explicitly, the
constants in this section depend on~$n$,~$\ell$ and our particular
choice of ${a_1,\ldots, a_{n-1}}$. In particular, {\it
sufficiently large} means {\it of absolute value  exceeding a
positive constant depending on~$n$,~$\ell$ and the choice of
${a_1,\ldots, a_{n-1}}$}.

One immediately verifies that ${f(0, 0) = 0}$.  Since ${f(0, x)}$
is a separable polynomial, it has a simple root at~$0$. Hence,
${\frac{\partial f}{\partial x}(0,0)\ne 0}$, and ${a_1\cdots
a_{n-1}\frac {\partial f}{\partial x}(0,0)}$ is a non-zero
integer.

Put
$$
a(t) := (-1)^{n-1}\frac{t^\ell-1}{a_1\cdots a_{n-1}}.
$$
Then
${f(t,x)=\left(x-a_1\right)\cdots\left(x-a_{n-1}\right)\left(x-
a(t)\right)-1}$. If ${\tau\in \Z}$ satisfies
\begin{equation}
\label{eaaa}
\tau\equiv 1 \mod a_1\cdots a_{n-1},
\end{equation}
then ${a(\tau) \in \Z}$ and
${f(\tau, x)\in \Z[x]}$. Moreover, for sufficiently large~$\tau$,
this polynomial gives rise to the ABC-field~$K_{a(\tau)}$, as
defined in Proposition~\ref{pabc}.

If the polynomial ${f(\tau, x)}$ has symmetric Galois group
over~$\Q$, then the field~$K_{a(\tau)}$ is primitive. Hence, the
following statement is a direct consequence of
Proposition~\ref{pabc}:\ref{iindex}.

\begin{proposition}
\label{pind} There exists a positive integer~$N$ (depending
on~$n$,~$\ell$ and the choice of ${a_1, \ldots, a_{n-1}}$) with
the following property. If the polynomial ${f(\tau, x)}$ has
symmetric Galois group over~$\Q$,  and the discriminant of the
field~$K_{a(\tau)}$ exceeds~$|\tau|^{2/(n+2)}$, then the index of
the group of ABC-units in the group of all units does not
exceed~$N$. \qed
\end{proposition}

\begin{definition}
\label{dsui}
An integer~$\tau$ is called {\it suitable} if it
satisfies the following conditions.

\begin{enumerate}

\item
\label{icong} We have~(\ref{eaaa}) and ${\gcd\left (\tau,
a_1\cdots a_{n-1}\frac {\partial f}{\partial x}(0,0)\right)=1}$.

\item
\label{isym}
The Galois group of the polynomial $f(\tau, x)$
(over~$\Q$) is symmetric.

\item
\label{idisc} The discriminant of~$K_{a(\tau)}$
exceeds~$|\tau|^{2/(n+2)}$.

\item
\label{ispecthin} Let~$\xi$ be the root of
${f(t,x)=\left(x-a_1\right)\cdots\left(x-a_{n-1}\right)\left(x-
a(\tau)\right)-1}$, as defined in
Proposition~\ref{pabc}:\ref{iroots}. Then, for all integers ${\nu,
\nu_1, \ldots, \nu_{n-1}} $, satisfying
\begin{equation}
\label{enus} 1\le \nu \le N, \quad 0\le \nu_1, \ldots, \nu_{n-1} <
\nu\ell,
\end{equation}
where~$N$ is defined in Proposition~\ref{pind}, and for every
determination of
$$
\zeta =\left(\xi\left(\xi-a_1\right)^{\nu_1/\nu}\cdots
\left(\xi-a_{n-1}\right)^{\nu_{n-1}/\nu}\right)^{1/\ell},
$$
we have ${[K(\zeta)\colon K]\ge \ell}$.

\end{enumerate}

\end{definition}

\begin{theorem}
\label{tsu} Put ${\mu =\frac1{2\ell(n-1)}}$. There exists a
positive constant~$c$ (depending on~$n$,~$\ell$ and the choice of
${a_1, \ldots, a_{n-1}}$) with the following property: for a large
positive real number~$X$ there exist at least~$cX^\mu$ suitable
integers~$\tau$
 which give rise to pairwise
non-isomorphic fields~$K_{a(\tau)}$ of discriminants not
exceeding~$X$.
\end{theorem}

\proof By Proposition~\ref{pabc}:\ref{idis}, there exists a
constant~$c_1$ with the following property: for any ${\tau\in
\Z}$, satisfying ${|\tau|\le c_1X^\mu}$, the discriminant of the
field~$K_{a(\tau)}$ does not exceed~$X$. Put ${T:=c_1X^\mu}$. Then
at least ${c_2T}$ integers~$\tau$ satisfy
\begin{equation}
\label{et1} \left(\frac{T^\ell+2}2\right)^{1/\ell} \le |\tau|\le T
\end{equation}
and condition~\ref{icong} of Definition~\ref{dsui}.

Integers~$\tau$ not satisfying  conditions~\ref{isym}
and~\ref{ispecthin} of Definition~\ref{dsui} form a thin set
(see Theorem~\ref{tgal} and
Proposition~\ref{pkum}). By Theorem~\ref{thit}, the number of
such~$\tau$ with ${|\tau|\le T}$ is ${O\left(\sqrt T\log
T\right)}$. Inequality~(\ref{edib}) implies that at most
${O\left(\sqrt T\right)}$ integers~$\tau$  with ${|\tau|\le T}$ do
not satisfy condition~\ref{idisc}. It follows that at least
${c_3T}$ suitable integers~$\tau$ satisfy~(\ref{et1}).

Finally,~(\ref{et1}) implies that ${A\le |a(\tau)|\le 2A}$, where
$$
A=\frac{T^\ell+1}{2\left|a_1\cdots a_{n-1}\right|}.
$$
By Proposition~\ref{pspr}, each~$K_{a(\tau)}$ may occur at most
${n(n-1)(n-2)}$ times. Hence the theorem is proved with
${c=c_3/(n(n-1)(n-2))}$. \qed

\section{The ABC-field corresponding to a suitable integer}
\label{sabcsi} \addtoc We are ready to complete the proof of
Theorem~\ref{tm}. In view of Theorem~\ref{tsu}, it remains to
prove the following.

\begin{proposition}
\label{peo} Let~$\tau$ be a suitable integer. Then the class group
of the field~$K_{a(\tau)}$ has an element of exact order~$\ell$.
\end{proposition}

\proof Since the suitable integer~$\tau$ is fixed, we may omit the
index and write ${K=K_{a(\tau)}}$. Since ${f(0,0)=0}$, the
polynomial ${f(0, x)}$ is divisible by~$x$. Put ${g(x) = a_1\cdots
a_{n-1} f(0,x)/x}$. Then ${g(x) \in \Z[x]}$ and ${g(0)= a_1\cdots
a_{n-1} \frac {\partial f}{\partial x}(0,0)}$. In particular,
${g(0) \ne 0}$ and
\begin{equation}
\label{egcd}
\gcd (g(0), \tau)=1.
\end{equation}

Rewrite the equality ${f(\tau, \xi)=0}$ as
$$
\xi g(\xi) = (-1)^{n-1}\left(\xi-a_1\right)\cdots\left(\xi -
a_{n-1}\right)\tau^\ell.
$$
 Since
${\xi-a_1, \ldots, \xi-a_{n-1}}$
 are units, this implies the
following equality for principal ideals:
${(\xi)(g(\xi))=(\tau)^\ell}$. Relation~(\ref{egcd}) implies
that~$\xi$ and~$g(\xi)$ are coprime. Hence, each of the principal
ideals~$(\xi)$ and~$(g(\xi))$ is an $\ell$-th power of an ideal
of~$K$.

Let~$\gera$ be the ideal of~$K$ such that ${\gera^\ell=(\xi)}$.
The order~$\lambda$ of the class of~$\gera$ in the class group
divides~$\ell$, and {\bf we wish to prove that ${\lambda=\ell}$.}

The ideal~$\gera^\lambda$ is principal. Fix ${\alpha \in K}$ such
that ${\gera^\lambda=(\alpha)}$ and let~$\zeta$ be some
determination of ${\alpha^{1/\lambda}}$. Then
\begin{equation}
\label{elamb}
[K(\zeta)\colon K]\le \lambda\le \ell.
\end{equation}
Let~$\nu$ be the index of ABC-units in the group of all Dirichlet
units of the field~$K$. Then any unit of~$K$ can be presented as
(a suitable determination of)
${\left(\xi-a_1\right)^{\nu_1/\nu}
\cdots\left(\xi-a_{n-1}\right)^{\nu_{n-1}/\nu}}$,
where ${\nu_1, \ldots, \nu_{n-1} \in \Z}$. In particular, since
${\zeta^{\ell} \in K}$ and ${\left(\zeta^\ell\right)= (\xi)}$, we
have
\begin{equation}
\label{ezetal} \zeta^\ell =
\xi\left(\xi-a_1\right)^{\nu_1/\nu}\cdots\left(\xi-a_{n-1}
\right)^{\nu_{n-1}/\nu}.
\end{equation}
Multiplying~$\zeta$ by a suitable ABC-unit, we may assume that the
integers ${\nu_1, \ldots, \nu_{n-1}}$ in~(\ref{ezetal}) satisfy
${0\le \nu_1, \ldots, \nu_{n-1} < \nu\ell}$. Also, since the
Galois group of the polynomial $f(\tau, x)$ is symmetric,
Proposition~\ref{pind} implies that ${\nu \le N}$.

Thus,~$\zeta$ satisfies~(\ref{ezetal}), where the integers ${\nu,
\nu_1, \ldots, \nu_{n-1}}$ satisfy~(\ref{enus}). Hence,
${[K(\zeta)\colon K]\ge \ell}$. Together with~(\ref{elamb}), this
implies ${\lambda=\ell}$, as wanted. This completes the proof of
Proposition~\ref{peo} and of Theorem~\ref{tm}. \qed

\section{Final remarks}
\addtoc
\begin{enumerate}
\item
Though we do not specifically consider the effectivity aspect  in
this note, one may check that our argument effectively bounds the
constants $X_0$ and~$c$ from Theorem~\ref{tm} in terms of~$\ell$
and~$n$.

\item
The estimate ${\left|\calf_n^{(\ell)}\right|\gg X^\mu}$ can be,
probably, slightly refined by letting the parameters ${a_1,
\ldots, a_{n-1}}$ vary.

\item
Theorem~\ref{tm} can be refined to count number fields with a
given number of real and complex embeddings. For this purpose,
one should simply replace our totally real ABC-fields by fields
with~$r$ real and~$2s$ complex embeddings, defined by polynomials
of the type
$$
(x-a_1)\cdots (x-a_r)\left(x^2+b_1x+c_1\right)\cdots
\left(x^2+b_sx+c_s\right)\pm1
$$
with ${b_j^2-4c_j<0}$.

\end{enumerate}

Some of these points  will be addressed in the forthcoming Ph.D.
thesis of S.~Hern\'andez.


\bigskip
{\small
\noindent
\begin{tabular}{lll}
Yuri F. Bilu & \hspace*{2cm}  &  Florian Luca \\
A2X & & Mathematical Institute \\
Universit\'e de Bordeaux 1 & & UNAM \\
351 Cours de la Lib\'eration & & Ap. Postal 61-3 (Xangari) \\
33405 Talence & & CP 58 089, Morelia, Michoac\'an\\
FRANCE & & MEXICO \\
yuri@math.u-bordeaux.fr & & fluca@matmor.unam.mx\\
\end{tabular}

}

}

\end{document}